\documentclass{article}

\usepackage{amsmath}
\usepackage{amsfonts}
\usepackage{amssymb}
\usepackage{mathrsfs}
\usepackage{graphicx}

\numberwithin{equation}{section}

\def \tab  {\hspace*{0.5cm}}
\def \ltab  {\hspace*{-0.5cm}}
\def \eud   {\, \mathrm d}

\def \im1  {\, \dot{\iota}\, }

\def \shiftright  {\hspace*{2.0cm}}
\def \shiftleft  {\hspace*{-2.0cm}}

\newcommand {\vect}[1]{ {\mathbf {#1}} }
\newcommand {\field}[1]{ {\mathbf {#1}} }

\newcommand {\Error}  { {\mathscr {\LARGE E}} }
\newcommand {\error}  { {\Large { \varepsilon}  } }
\newcommand {\SquareIntegrable}    {{\mathscr {L}}^2}
\newcommand {\reals} {{\mathbb {R}}}

\newcommand {\e}{\textit{\Large e}}
\newcommand {\Continuous}  { {\mathscr {\LARGE C}} }

\begin{document}

\title{Dynamic Programming Method for Best Piecewise Linear Approximation for Vector Field of Nonlinear Boundary Value Problems on the Interval $[0, \, 1]$}

\author{\bf{Duggirala Meher Krishna}\\
{\small{Gayatri Vidya Parishad College of Engineering (Autonomous)}} \\
{\small{Madhurawada, VISAKHAPATNAM -- 530 048, Andhra Pradesh, India}} \\
 {\small{E-mail ~: \tab duggiralameherkrishna@gmail.com}}\\
 \\
 and \\
 \\
\bf{Duggirala Ravi}\\
{\small{Gayatri Vidya Parishad College of Engineering (Autonomous)}} \\
{\small{Madhurawada, VISAKHAPATNAM -- 530 048, Andhra Pradesh, India}} \\
\shiftleft \tab {\small{E-mail ~: \tab ravi@gvpce.ac.in; \tab duggirala.ravi@yahoo.com};} \\
{\small{\shiftright  duggirala.ravi@rediffmail.com; \tab drdravi2000@yahoo.com}} 
}

\date{}

\maketitle

\begin{abstract}
An important problem that arises in many engineering applications is the boundary value problem for ordinary differential equations. There have been many computational methods proposed for dealing with this problem. The convergence of the iterative schemes to a true solution, when one such exists, and their numerical stability are the central issues discussed in the literature. In this paper, we discuss a method for approximating the vector field, maintaining the boundary conditions and numerical stability. If a true solution exists, finer discretization of the solution space converges to one such.  
\end{abstract}

\section{Introduction}

During the last few decades there has been a remarkable growth of interest in problems associated with solving linear and nonlinear ordinary differential equations satisfying boundary conditions. For many of the nonlinear boundary value problems that occur in engineering and applied sciences, it is difficult to obtain a solution analytically. For a nonlinear boundary value problem, the difficulty lies in establishing the existence of a solution mathematically, though in some cases multiple solutions exist. Approximation of the solution space of a given differential equation has gained importance as it speeds up or helps in solving the problem efficiently. These approximation methods can be put into two classes: (i) those in which a solution is approximated numerically at a number of discrete points of the domain, and (ii) those in which a solution is approximated by a finite number of terms of a sequence of functions. The approach in (ii) is called a weighted residual method.  In most of the numerical methods described in the literature, we may have to add a sufficient number of some undetermined variables with implicitly assumed conditions at one end of the domain, and adjust the additional variables until the required conditions are satisfied at the other end to obtain the solution of the boundary value problem [\cite{Daniel1978} -- \cite{Fox1990},  \cite{Keller1987},  \cite{Osbourne1969} --  \cite{RS1972}], and further approximate the derivatives of the dependent variables with forward, backward or central difference operators defined at the grid points \cite{NS1972}. In the first part of a numerical scheme as just mentioned, the convergence may be very slow, and in the second part, convergence and stability of the particular difference scheme may depend on the selection of the approximations used for the derivatives involved in the differential equation and the boundary conditions. It may also be the case that the chosen difference method is not numerically stable, resulting in chaos phenomenon creeping into the iterative schemes, at places where matrix inversions are utilized, without the solver being explicitly aware of its entry. In some cases, suitable regularization and relaxation conditions, involving more variables, may have to be added to the constraints formulated in the previous steps. In the weighted residual methods, difference equations are generated using approximation methods with piecewise polynomial solutions \cite{NS1972}. 

Among the most popular and successful techniques for solving boundary value problems with nonlinearities is Galerkin procedure. In this approach, the solution of the ordinary differential equation is expressed as a linear combination of certain basis functions, and the coefficients of the basis functions are determined by requiring that the residual be orthogonal to each of the basis functions. The difficulty lies in the selection of basis functions to obtain the desired solution, that can take care of the boundary conditions simultaneously. In recent times, the concept of piecewise linear approximation of the differential equation gained momentum [\cite{GR1996}, \cite{Keller1987}, \cite{NS1972} -- \cite{Osbourne1969}] The two point boundary value problems are approximated by piecewise linear ones which have analytical solutions and reduced to finding the slope of the solution at the left boundary, so that the boundary conditions at the right end of the interval are satisfied. This approach results in a complex system of non-linear algebraic equations. Some more recent and highly efficient algorithms \cite{SUNDIALS}, for solving these complex systems of differential and algebraic equations (DAE) can be used for computational purpose. 

The motivation for the present work is the necessity of a new method that is applicable for most or all general continuous vector fields and general boundary conditions. The objective of our study is to find efficiently a solution, if exists, by an algorithm, such that the algorithm is able to detect and report when there is no solution, if an error term does not fall below a threshold, despite using various approximation schemes with several basis functions. We propose a new method based on dynamic programming for solving boundary value problems in one variable. The dynamic programming based formulation is adapted for obtaining an optimal approximation for the vector field for genera of two-point boundary value problem, which is usually formulated as an optimal control problem in the literature \cite{RK1965}. For improving an initial approximation, repeated application of the dynamic programming algorithm with refined discretization of the parameter space can be used. A modified Newton-Raphson method for improving an approximate solution along with updating the initial value is also discussed. These aspects are newly introduced in this work.  Possible extensions to formulation of solution methods for optimal control problem and for boundary value problems of partial differential equations defined on compact simply connected convex domains are briefly described.

\section{\label{Section-DPM}Best Piecewise Linear Approximations for Vector Fields of Two-point Boundary Value Problems}

In this section, we consider the boundary value problem
\begin{equation}
\left.
\begin{array}{cl}
   \vect{x}'(t) = \field{f}(\vect{x}(t),\, t), & \textrm{ for }   0 < t < 1,\\
   \beta(\vect{x}(0), \, \vect{x}(1)) = 0,  &  
\end{array}
\right\} \label{BVP-Form1}
\end{equation}
where $\field{f} \in \Continuous(\reals^n \times [0,\, 1]; \, \reals^n)$,
$\beta \in \Continuous(\reals^n \times \reals^n ;\,  \reals)$, 
and $ \vect{x}'(t) = \frac{\eud\vect{x}}{\eud t} $.
The vectors are written as columns of appropriate
dimensions. The objective of obtaining a piecewise
linear approximation is stated as follows:
\begin{enumerate}
\item Let, for some fixed integer $N \geq 2$, 
   $0 = t_0 < t_1 < \ldots < t_N = 1$ be
fixed node points of the time interval $[0, \, 1]$.

\item Then, we obtain optimal values for the
parameters of form $\pi_i = (A_i(t), \, \vect{b}_i(t), \, \theta_i)$,
 $0 \leq i \leq N-1$, where $A_i$  and  $\vect{b}_i$
are $n \times n$ and $n\times 1$ matrices
consisting of undetermined linear combination
of some fixed basis functions defined on
$[t_i, \, t_{i+1}]$, and $\theta_i$ is $n \times 1$
column vector of undetermined constants,  
such that the curve $\vect{y}(t)$ satisfying 
\begin{equation}
\left.
\begin{array}{l}
\vect{y}'(t) = A_i (t) \vect{y}(t) + \vect{b}_i(t), 
~ \textrm{ for } t_i < t < t_{i+1} \nonumber \\
 \textrm{ with the initial condition }
\vect{y}(t_i) = \theta_i
\end{array}
\right\} \textrm{ for } 0 \leq i \leq N-1
  \label{ArcPiece-1}
\end{equation}
 subject to the boundary and continuity conditions
\vspace*{-0.2cm}
\begin{equation}
\left.
\ltab
\begin{array}{cl}
  \beta(\theta_0, \, \theta_N) ~ = ~ 0,  & \textrm{ where } \theta_N = \vect{y}(1^{-}) \\
 \vect{y}(t_{i+1}^-) ~ = ~  \theta_{i+1},  &  \textrm{ for } 0 \leq i \leq N-1, \textrm{ and } \\
\ltab  A_i(t_{i+1}) \theta_{i+1} + \vect{b}_i(t_{i+1}) &
 \ltab~ =  A_{i+1}(t_{i+1}) \theta_{i+1} +  \vect{b}_{i+1}(t_{i+1}),\\
&   \textrm{ for } 0 \leq i \leq N-2
\end{array}
\right \} \label{BCC-General}
 \end{equation}
minimizes the error defined by
\begin{equation}
\Error = \sum_{i = 0}^{ N-1} \int_{t_i}^{t_{i+1}} \mathcal{G}
 {\bigg (} A_i(t) \vect{y}(t) + \vect{b}_i(t) - \field{f}(\vect{y}(t),\, t) {\bigg )} \eud t
    \label{ErrorFunctional-1}
\end{equation}
where the integrand $\mathcal{G}$ can be chosen to be
any nonnegative continuous function defined on $\reals^n$
and equal to 0 at the origin. 
\end{enumerate}

We first consider a special case of (\ref{BVP-Form1}),
in which the boundary conditions are variables
seperable (i.e. of the form $\beta_0(\vect{x}(0)) = 0$
and $\beta_1(\vect{x}(1)) = 0$, where
$\beta_0, \, \beta_1 \in \Continuous(\reals^n;\, \reals)$).
Later, we shall indicate the modifications
needed to extend the method to the general case of
boundary conditions described in (\ref{BVP-Form1}).

\subsection{Piecewise Linear Approximation for the Special Case}

In this section, we consider the boundary value problem
\begin{eqnarray}
&&    \vect{x}'(t) = \field{f}(\vect{x}(t),\, t),
 \tab \textrm{ for }   0 < t < 1,\nonumber\\
&&   \beta_0(\vect{x}(0)) = \beta_1(\vect{x}(1))  =  0,   
 \label{BVP-SpecialCase}
\end{eqnarray}
where $\field{f} \in \Continuous(\reals^n\times[0,\, 1] ; \, \reals^n)$, and
$\beta_0, \, \beta_1 \in \Continuous(\reals^n \times \reals^n ;\,  \reals)$,
where $\Continuous(X ; \, Y)$ is the space of continuous functions from $X$ into $Y$
and $\reals$ is the set of real numbers. The objective being described in the sequel
is to find optimal values for the parameters
$\pi_i = (A_i, \, B_i, \, \theta_i)$, $0 \leq i \leq N-1$, for which
the error functional in (\ref{ErrorFunctional-1}) is minimized by the solution
satisfying (\ref{ArcPiece-1}) subject to the boundary and continuity conditions
\vspace*{-0.2cm}
\begin{equation}
\left.
\begin{array}{cl}
  \beta_0(\theta_0) \, = \, 0, ~ \beta_1(\theta_N) \, = \, 0, & \textrm{ where } \theta_N = \vect{y}(1) \\
 \vect{y}(t_{i+1}^-) ~ = ~  \theta_{i+1},  &  \ltab \textrm{ for } 0 \leq i \leq N-1, \textrm{ and } \\
\shiftleft A_i(t_{i+1}) \theta_{i+1} + \vect{b}_i(t_{i+1}) &
\ltab\ltab \ltab~ = ~   A_{i+1}(t_{i+1}) \theta_{i+1} +  \vect{b}_{i+1}(t_{i+1}),\\
&\textrm{ for } 0 \leq i \leq N-2
\end{array}
\right \} \label{BCC-SpecialCase}
 \end{equation}
 The solution for (\ref{ArcPiece-1})
in the interval $[t_i, \, t_{i+1})$,
 $0 \leq i \leq N-1$,  is
\begin{equation}
 \vect{y}(t)   =   \e ^{{\bigg[}\int_{t_i}^{t}A_i(s)\, \eud s{\bigg]}} \, \theta_i  +
	 \int_{t_i}^{t} \e ^{{\bigg[}\int_{s}^{t}A_i(v)\, \eud v{\bigg]}} \, \vect{b}_i(s) \, \eud s, 
		\label{Solution-ArcPiece-1}
\end{equation}
which is expressible explicitly in terms
of the parameter values. Substituting
the expression for $\vect{y}(t)$ in
(\ref{ErrorFunctional-1}), the value
of $\Error$ can be found for any
prescribed values of the parameters.
If the parameter space is discretized,
the values of $\Error$ can be
tabulated for various values of
the parameters. However,
using the additive property of
the error functional $\Error$,
it is possible to formulate a
dynamic programming method
for obtaining the tables quickly
and efficiently. For this purpose,
we define a $k$-step error
functional $\Error_k$, for
$0 \leq k \leq N$, by the
relation $\Error_0 = 0$ and
for $k = 1, 2, \ldots, N$,
\begin{eqnarray}
&& \error_k~ = ~  \int_{t_{k-1}}^{t_{k}} \mathcal{G}
{\bigg (} A_{k-1}(t) \vect{y}(t) + \vect{b}_{k-1}(t) -
 \field{f}(\vect{y}(t),\, t) {\bigg )} \eud t     \label{error-k}\\
&& \Error_k ~ = ~ \Error_{k-1} +\error_k
    \label{ErrorFunctional-2}
\end{eqnarray}
Thus $\Error_N$ in (\ref{ErrorFunctional-2})
is actually $\Error$ in (\ref{ErrorFunctional-1}).
The functional $\Error_k$ depends only on the
parameters $(\pi_0, \pi_1, \ldots, \pi_{k-1})$.
To make this dependency explicit, we sometimes
write $\Error_k(\pi_0, \pi_1, \ldots, \pi_{k-1})$ for 
$\Error_k$; similarly, $\error_k$ is sometimes
written as $\error_k(\pi_{k-1})$.  Now let
for $1 \leq k \leq N$,  $\Omega_k$ be the
set of parameters $(\pi_{k-1},\, \pi_k)$
satisfying
\begin{equation}
\left.
\begin{array}{c}
\beta_0(\theta_0) \, = \, 0, ~  \textrm{  if } k = 1 \\
  \beta_N(\theta_N) \, = \, 0, ~ \textrm{ if } k = N \\
  \vect{y}(t_k^-)  ~ = ~ \theta_k,  ~ \textrm{ and } \\
\ltab \ltab  A_{k-1}(t_{k}) \theta_{k} + \vect{b}_{k-1}(t_{k})  \, = \, A_{k}(t_{k}) \theta_{k} +  \vect{b}_{k}(t_{k}),
 ~ \textrm{ if } k < N \\
\end{array}
\right \}  
\label{Omega-k-SpecialCase}
\end{equation}
where if $k = N$, the parameter $\pi_N$ is
interpreted as $\theta_N$, and $\vect{y}(t)$
in the interval $[t_{k-1},\, t_k)$ is given by
(\ref{Solution-ArcPiece-1}) with $i = k-1$.
Let
 \begin{equation}
\Gamma_k = \{(\pi_0, \pi_1, \ldots, \pi_k)\, | \;
   (\pi_{i-1},\, \pi_i) \in \Omega_i , \; \textrm{ for } 1\leq i \leq k \} 
\label{Gamma-k-SpecialCase}
\end{equation}
Now we define a finite sequence of functions
on the parameter space as follows: let
$S_0(\pi_0) = 0$, and for $1 \leq k \leq N$,   
\begin{equation}
S_k(\pi_k) \tab  =  \min_{ \scriptsize{ \begin{array}{c}
(\pi_0, \ldots, \pi_{k-1})\\
~\textrm{\footnotesize such that}\\
(\pi_0, \ldots, \pi_{k-1}, \pi_k) \,\in\, \Gamma_k
\end{array}} } \ltab   \Error_k(\pi_0, \ldots, \pi_{k-1}) 
\label{MinimalFunction-Sk-SpecialCase}
\end{equation}
Then, observing that for each fixed $\pi_k$,
\begin{small}
\[
 \min_{ (\pi_0, \ldots, \pi_{k-1}, \pi_k) \in \Gamma_k}
                     \Error_k(\pi_0, \ldots, \pi_{k-1}) =
 \min_{  \scriptstyle{ ^{ (\pi_0,  \ldots, \pi_{k-1}) \in \Gamma_{k-1}}_{
               \textrm{\footnotesize and }  (\pi_{k-1},\,\pi_k) \in \Omega_k}} } \{
         \Error_{k-1}(\pi_0,  \ldots, \pi_{k-2}) + \error_k(\pi_{k-1}) \}
\]
\end{small}
 we can conclude that
\begin{equation}
S_k(\pi_k)\tab  = \min_{ \scriptstyle{ ^{\pi_{k-1} ~\textrm{\footnotesize such that} }_{
          ~ (\pi_{k-1}, \, \pi_k) \in \Omega_k} } }  \{ S_{k-1}(\pi_{k-1}) + \error_k(\pi_{k-1})\}
\label{RecurrenceRelation-SpecialCase}
\end{equation}
The value of $\min_{\pi_N} S_N(\pi_N) =
\min_{\theta_N} S_N(\theta_N)$ gives
the minimum value of $\Error_N$ over the
parameter space constrained by the boundary
and continuity conditions (\ref{BCC-SpecialCase}).

Now we describe a tabulation procedure for
computation of optimal parameters as follows:
\begin{enumerate}
\item Initially compute $S_1(\pi_1) = \min_{(\pi_0,\, \pi_1) \in \Omega_1} \error_1(\pi_0)$,
and set 
\[
\pi_0^{\ast}(\pi_1) = \{\pi_0\, | \; (\pi_0,\, \pi_1) \in \Omega_1 \textrm{ and } \error(\pi_0) = S_1(\pi_1)\}
\]

\item Now for $2 \leq k \leq N$, compute
$S_k(\pi_k) = \min_{(\pi_{k-1},\, \pi_k) \in \Omega_k} S_{k-1}(\pi_{k-1})+\error_k(\pi_{k-1})$,
and set
\[
\pi_{k-1}^{\ast}(\pi_k) =
\{\pi_{k-1}\,|\; (\pi_{k-1},\, \pi_k) \in \Omega_k \textrm{ and } S_k(\pi_k) =
S_{k-1}(\pi_{k-1})+\error_k(\pi_{k-1})\}
\]

\end{enumerate}

At the end of the procedure, we obtain a table for $S_N(\pi_N)$
for various values of $\pi_N$ (which is $\theta_N$). The optimal
parameters are found by backtracking as follows:
\begin{enumerate}
\item First find $\pi_N^{\ast} = \textrm{argmin}_{\pi_N} S_N(\pi_N)$.

\item Then for $k = N-1,\,  \ldots, \, 1,\, 0$, choose $\pi_k^{\ast}$
as an element of $\pi_k^{\ast}(\pi_{k+1}^{\ast})$.

\end{enumerate}
The tuple $(\pi_0^{\ast},\, \pi_1^{\ast},\, \ldots,\, \pi_N^{\ast})$
thus obtained are optimal parameters with respect to
the error functional (\ref{ErrorFunctional-1}) subject
to the boundary and continuity conditions in (\ref{BCC-SpecialCase}).

\subsection{Piecewise Linear Approximation for the General Case}

Now we describe a generalization of the method of the
previous section to the general case of boundary conditions
in (\ref{BCC-General}). Let for $1 \leq k \leq N$, 
$\Omega_k$ be the set of parameters $(\pi_{k-1},\, \pi_k)$
satisfying
\begin{equation}
\left.
\begin{array}{c}
  \exists{\xi_N}: \,\beta(\theta_0, \xi_N) \, = \, 0, ~  \textrm{  if } k = 1 \\
\exists{\xi_0}: \, \beta(\xi_0, \theta_N) \, = \, 0, ~  \textrm{  if } k = N \\
  \vect{y}(t_k^-)  ~ = ~ \theta_k,  ~ \textrm{ and } \\
\ltab \ltab  A_{k-1}(t_{k}) \theta_{k} + \vect{b}_{k-1}(t_{k})  \, = \, A_{k}(t_{k}) \theta_{k} +  \vect{b}_{k}(t_{k}),
 ~ \textrm{ if } k < N \\
\end{array}
\right \}  
\label{Omega-k-GeneralCase}
\end{equation}
Let $\Gamma_k$ be as in (\ref{Gamma-k-SpecialCase}),
but with $\Omega_{k}$ as in (\ref{Omega-k-GeneralCase}).
Now the parameter space constrained by (\ref{BCC-General})
is given by
$\Gamma_N \bigcap {\mathcal V}_N$,  where
${\mathcal V}_N = \{(\pi_0,\, \pi_1, \, \ldots, \, \pi_N)\, |
 \;\beta(\theta_0,\, \theta_N) = 0\}$.
However, $S_k$ is now defined on $(\pi_0, \, \pi_k)$:
\begin{equation}
S_k(\pi_0,\, \pi_k) \tab = \min_{\scriptstyle { ^{ (\pi_1, \ldots, \pi_{k-1}) ~
       \textrm{\footnotesize such that } }_{
      ~ (\pi_0,\, \pi_1,\,\ldots,\,\pi_{k-1},\, \pi_k) \in \Gamma_k } } }
                     \Error_k(\pi_0, \pi_1, \ldots, \pi_{k-1})
\label{MinimalFunction-Sk-GeneralCase}
\end{equation}
so that
\begin{equation}
S_k(\pi_0,\, \pi_k) \tab = \min_{ \scriptstyle{ ^{\pi_{k-1} ~\textrm{\footnotesize such that} }_{
          ~ (\pi_{k-1}, \, \pi_k)\, \in \, \Omega_k} } }  \{ S_{k-1}(\pi_0,\, \pi_{k-1}) + \error_k(\pi_{k-1})\}
\label{RecurrenceRelation-GeneralCase}
\end{equation}
for $2 \leq k \leq N$. Now 
\begin{equation}
\ltab 
\min_{ \scriptsize{\begin{array}{c}
 (\pi_0, \, \pi_N)\\
 \textrm{ such that }\\ 
 \beta(\theta_0,\,\theta_N) = 0\\
 \end{array}}}  \ltab S_N(\pi_0, \pi_N) ~~  = 
\ltab \ltab \min_{\scriptsize {\begin{array}{c}
(\pi_1,\, \ldots, ,\,\pi_{N-1})\\
       \textrm{ such that }\\
      (\pi_0,\, \pi_1,\,\ldots,\, \pi_N) \, \in \,  {\mathcal V}_N \bigcap \Gamma_N\\
\end{array} } } \ltab \ltab \Error_N(\pi_0, \pi_1, \ldots, \pi_{N-1})
                \label{MinimumEquivalence-GeneralCase}
\end{equation}
The computation of the forward tables is
as follows:
\begin{enumerate}
\item Initially compute $S_1(\pi_0,\,\pi_1)$ for various values of $(\pi_0, \,\pi_1) \in \Omega_1$.

\item Now for $2 \leq k \leq N$, compute
$S_k(\pi_0, \, \pi_k)$ from
(\ref{RecurrenceRelation-GeneralCase}),
and set
\begin{eqnarray*}
\pi_{k-1}^{\ast}(\pi_0,\, \pi_k)  & = &
\{\pi_{k-1}\,|\; (\pi_{k-1},\, \pi_k) \in \Omega_k \textrm{ and }\\
& & \tab \tab S_k(\pi_0,\, \pi_k) = 
S_{k-1}(\pi_0,\, \pi_{k-1})+\error_k(\pi_{k-1})\}
\end{eqnarray*}

\end{enumerate}
From the table thus computed, the optimal values
of the parameters are extracted by backtracking
as follows:
\begin{enumerate}
\item First find $(\pi_0^{\ast},\,\pi_N^{\ast})$ from
\[
 \begin{array}{rccl}
 (\pi_0^{\ast},\,\pi_N^{\ast})  & = & \textrm{argmin} &  S_N(\pi_0,\, \pi_N)\\
&& \scriptstyle{(\pi_0,\,\pi_N)} ~\textrm{\footnotesize such that}&\\
&&\scriptstyle{B(\theta_0,\, \theta_N) \, = \, 0}& \\
\end{array}
\]

\item Then for $k = N-1,\,  \ldots, \,2,\, 1$, choose $\pi_k^{\ast}$
as an element of $\pi_k^*(\pi_0^{\ast},\, \pi_{k+1}^{\ast})$.

\end{enumerate}
The tuple $(\pi_0^{\ast},\, \pi_1^{\ast},\, \ldots,\, \pi_N^{\ast})$
thus obtained minimizes the error functional (\ref{ErrorFunctional-1})
subject to the boundary and continuity conditions in (\ref{BCC-General}).

\section{Improving the Initial Approximation}

In this section, we describe a method similar
to Newton-Raphson method for improving
an initial approximation for a BVP or
in particular, for an initial value problem.
The correction of the approximation consists
of two parts: in the first part, a correction term
with 0 initial value for improving the solution in
the interior is obtained, and the second part
finds optimal correction term in the initial value
so that the updated solution satisfies the boundary
conditions more accurately. Alternatively, it is also
possible to achieve the same by either Picard's
successive approximation (with suitable correction
term in the initial value) or repeated application
of the dynamic programming method described
in the last section with finer discretization of the
parameter space restricted to a tube-like set
around the initial approximate solution. 

\subsection{\label{Section-Improvements}Improving an Approximate Solution with Intial Value Fixed}
We assume that the vector field  
$\field{f}$ is differentiable with
continuous derivatives of upto as
high an order (upto two) as necessary.
Let $\vect{x}_0(t)$ be an initial approximate
solution satisfying 
\begin{equation}
\vect{x}_0'(t) \approx \field{f}(\vect{x}_0(t),\, t), 
 ~ \textrm{ for } 0 < t < 1 \label{Approximation0}
\end{equation}
The objective is to formulate an efficient method
for improving the approximation in (\ref{Approximation0})
by successive iterations.  Let for $k = 0, \, 1, \, 2,\,  \ldots$, 
\begin{equation}
\left.
\begin{array}{c}
\vect{x}_k'(t) = \field{f}_k(t), 
 ~ \textrm{ for } 0 < t < 1\\
\textrm{with the initial condition } \vect{x}_k(0) = \theta_0
\end{array}
\right \}  \label{Equation-k}
\end{equation}
At step $k+1$, the correction term
$\vect{y}_k = \vect{x}_{k+1} - \vect{x}_k$
must be found such that
\begin{equation}
\left.
\begin{array}{c}
\vect{x}_{k+1}'(t) = \field{f}(\vect{x}_{k+1}(t),\, t)
 ~ \textrm{ for } 0 < t < 1\\
\textrm{with the initial condition } \vect{x}_{k+1}(0) = \theta_0
\end{array}
\right \} \label{ExactSolution}
\end{equation}
Subtracting (\ref{Equation-k})
from (\ref{ExactSolution}), we 
find that $\vect{y}_k$ must satisfy
\begin{equation}
\left.
\begin{array}{c}
\vect{y}_k'(t) = \field{f}(\vect{x}_{k+1}(t),\, t)-\field{f}_k(t), 
 ~ \textrm{ for } 0 < t < 1\\
\textrm{with the initial condition } \vect{y}_k(0) = 0
\end{array}
\right \}  \label{Correction-Term-Equation-k}
\end{equation}
Now using the Taylor series approximation
with respect to the first variable (i.e. $\vect{x}$)
upto the first order for $\field{f}(\vect{x}_{k+1}(t),\, t)$
we find
\begin{eqnarray}
\field{f}(\vect{x}_{k+1}(t),\, t)  & = &
\field{f}(\vect{x}_{k}(t)+\vect{y}_{k}(t),\, t) \nonumber \\
& \approx & \field{f}(\vect{x}_{k}(t),\, t) + A_k(t)\vect{y}_{k}(t) 
 \label{TaylorSeriesApproximation-f}\\
& & \shiftright \textrm {where} ~ A_k(t) =
{\bigg[} \partial_{\vect{x}}\field{f}(\vect{x},\, s){\bigg]}_{
\scriptstyle{ ^{s = t}_{\vect{x} = \vect{x}_k(t) } } }
\nonumber
\end{eqnarray}
Letting $\field{b}_k(t) = 
 \field{f}(\vect{x}_{k}(t),\, t) - \field{f}_k(t)$
and using (\ref{TaylorSeriesApproximation-f})
in (\ref{Correction-Term-Equation-k}),
an approximate correction term
is found by solving
\begin{equation}
\left.
\begin{array}{c}
\vect{y}_k'(t) = A_k(t)\vect{y}_{k}(t)  + \field{b}_k(t), 
 ~ \textrm{ for } 0 < t < 1\\
\textrm{with the initial condition } \vect{y}_k(0) = 0
\end{array}
\right \}  \label{Approximate-Correction-Term-Equation-k}
\end{equation}
 The solution of (\ref{Approximate-Correction-Term-Equation-k})
is given by
\begin{equation}
\vect{y}_k(t) =  \int_{0}^{t} \e ^{{\bigg[}\int_{s}^{t}A_k(v)\,\eud v{\bigg]}} 
   \, \vect{b}_k(s) \, \eud s   \label{Approximate-Correction-Term-k}
\end{equation}
The iteration converges fast (at almost quadratic rate) 
provided the initial approximate solution
(\ref{Approximation0}) is sufficiently close
to the exact solution.  If the matrix $A_k(t)$ is
uniformly boundedly invertible in the interval
$[0,\,1]$, then $\vect{y}_k(t)$ can also be
taken as
\begin{equation}
\vect{y}_k(t) = -A_k^{-1}(t) \vect{b}_k(t), \tab 0 \leq t \leq 1 
 \label{Actual-NRM-for-Correction-Term}
\end{equation}
which is the well-known method for solving for $\vect{x}$
from $\field{f}(\vect{x},t) = 0$ for each $0 \leq t \leq 1$.
This method also modifies the initial condition. If the
initial value is required to be updated independently,
then we have to find $\vect{y}_k(t)$ from 
(\ref{Approximate-Correction-Term-Equation-k}),
with solution given by (\ref{Approximate-Correction-Term-k}).
It may be observed when either of
(\ref{Approximate-Correction-Term-k})
and (\ref{Actual-NRM-for-Correction-Term})
converges,
$\lim_{k \rightarrow \infty} \vect{b}_k(t) = 0$,
almost everywhere (a.e.), for $0 \leq t \leq 1$, which
implies, if convergent, either iteration leads
to the final solution satisfying
$\vect{x}'(t) = \field{f}(\vect{x}(t),\, t)$
a.e., for $0 \leq t \leq 1$.

\subsection{\label{Section-Modified-Improvements}Improving the Initial Value}
Suppose that the initial value of approximate
correction term $\vect{y}_k$ is
$\eta_k$ (to be determined). Then the
correction term is given by
\begin{equation}
\vect{y}_k(t;\,\eta_k) = \e ^{{\bigg[}\int_{0}^{t}A_k(v)\,\eud v{\bigg]}} \eta_k +
      \int_{0}^{t} \e ^{{\bigg[}\int_{s}^{t}A_k(v)\,\eud v{\bigg]}} 
   \, \vect{b}_k(s) \, \eud s   \label{InitialValue-Correction-Term-k}
\end{equation}
where $A_k(t)$ is as in
(\ref{TaylorSeriesApproximation-f})
and $\vect{b}_k(t)$ in the
following line. Then the
parameter $\eta_k$ is found based
on an optimality criterion. The objective
can be formulated as follows:
\begin{equation}
\left.
\begin{array}{c}
\eta_k^{\ast} = \textrm{argmin}_{\scriptsize{_{\eta_k}}}
\mathscr{F}(\eta_k), \textrm{ where}\\
\mathscr{F}(\eta_k)  =  \int_0^1
\| \vect{y}'_k(t ;\,\eta_k)\|^2 \, \eud t
\end{array}
\right\} \label{ObjectiveFunction-Eta-k}
\end{equation}
If the initial approximation is sufficiently accurate
at $t = 0$, then $\eta_k$ must be small. Thus
we can expect that  the minimum in
(\ref{ObjectiveFunction-Eta-k}) is attained for
 $\eta_k \approx 0$. The update value can be
chosen to be proportional to the gradient of
$\mathscr{F}(\eta_k)$ at 0. Specifically, 
we can choose
$\eta_k^{\ast} = - h_k \cdot \nabla \mathscr{F}(0)$,
for a small positive number $h_k$, resulting
in an easy updation of the initial value.
The multiplier $h_k$ can be found by
binary search method over an interval
of the form $[0, \, M]$, for some sufficiently
large positive constant $M$.  Alternately
it is also possible to find $\eta_k^{\ast}$
by solving $\nabla {\mathscr F}(\eta) = 0$,
which gives the iterative formula $\eta_0 = 0$
and for $k = 1,\, 2,\, \ldots$,
\begin{equation}
\eta_k = \eta_{k-1} - H^{-1}(\eta_{k-1}) \nabla {\mathscr F}(\eta_{k-1})
   \label{NRM-for-Initial-Value-Update}
\end{equation}
where $H(\eta)$ is the Hessian matrix of
${\mathscr F}(\eta)$. The method is fast
and does not require a seperate search
for the multiplier constant as in the case
of gradient descent method.

\subsection{Combining Both: Improving an Initial Approximation by Successive Iterations}
In this section we briefly describe a gradient
descent method for imporving an initial solution
of the boundary value problem (\ref{BVP-Form1}).
For this purpose, we assume that $\beta$ is
continuously differentiable having continuous derivatives
upto Hessian. The update in approximation $\vect{x}_k(t)$
of the $k$-th iteration is given by
(\ref{InitialValue-Correction-Term-k})
with initial condition $\vect{y}_k(0) = \eta_k$,
and $\vect{y}'_k(t)$ satisfying
(\ref{Correction-Term-Equation-k}).
The initial value $\eta_k$ of $\vect{y}_k(t)$ 
is determined by solving the boundary condition.
Specifically let $\theta_0^{(k)} = \vect{x}_k(0)$
and $ \theta_N^{(k)} = \vect{x}_k(1)$. The initial
value  $\theta_0^{(k+1)} = \vect{x}_{k+1}(0)$
is  $\theta_0^{(k)}+\eta_k$, and the final value
$\theta_N^{(k+1)}$ is $\vect{x}_k(1)+\vect{y}_k(1)$.
The updated values $(\theta_0^{(k+1)},\,\theta_N^{(k+1)})$
must be found such that 
\begin{equation}
\beta(\theta_0^{(k+1)},\,\theta_N^{(k+1)}) ~  =  ~ 0
\label{Beta-k+1-A}
\end{equation}
Now
\begin{eqnarray}
&&\theta_0^{(k+1)}  ~ = ~ \theta_0^{(k)} + \eta_k
\, , \tab \textrm{and}   \label{Theta-0-k+1}\\
&&\ltab \ltab \theta_N^{(k+1)} ~ = ~ \theta_N^{(k)} + 
\e ^{{\bigg[}\int_{0}^{1}A_k(v)\,\eud v{\bigg]}} \eta_k +
      \int_{0}^{1} \e ^{{\bigg[}\int_{s}^{1}A_k(v)\,\eud v{\bigg]}} 
   \, \vect{b}_k(s) \, \eud s 
 \label{Theta-N-k+1}
\end{eqnarray}
Substituting the values of 
$\theta_0^{(k+1)}$ and $\theta_N^{(k+1)}$
from (\ref{Theta-0-k+1}) and (\ref{Theta-N-k+1})
into (\ref{Beta-k+1-A}), we find
\begin{equation}
 \beta{\bigg(}\theta_0^{(k)}+\eta_k,\,
\theta_N^{(k)} +  \e ^{{\bigg[}\int_{0}^{1}A_k(v)\,\eud v{\bigg]}} \eta_k +
      \int_{0}^{1} \e ^{{\bigg[}\int_{s}^{1}A_k(v)\,\eud v{\bigg]}} 
   \, \vect{b}_k(s) \, \eud s {\bigg)} = 0
\label{Beta-k+1-B}
\end{equation}
As (\ref{Beta-k+1-B}) is to be solved for
the vector $\eta_k$ from only one equation,
we propose first a gradient descent method
for minimization of $\beta^2$. Further,
if $\eta_k$ is small, we can evaluate the
gradient of $\beta$ in (\ref{Beta-k+1-B})
with respect to $\eta_k$ for $\eta_k = 0$.
Thus $\eta_k$ is chosen such that
\begin{equation}
\eta_k = - h_k \cdot \beta(\theta_0, \, \theta_N)\cdot {\bigg[}
\partial_{\theta_0} \beta(\theta_0,\, \theta_N)+
\e ^{{\bigg[}\int_{0}^{1}A_k(v)\,\eud v{\bigg]}}\cdot 
\partial_{\theta_N} \beta(\theta_0,\, \theta_N){\bigg]}
\label{Eta-k-BVP}
\end{equation}
where $\beta(\theta_0,\, \theta_N)$ and its partial derivatives
$\partial_{\theta_0} \beta(\theta_0,\, \theta_N)$
and  $\partial_{\theta_N} \beta(\theta_0,\, \theta_N)$
are evaluated for $\theta_0 = \theta_0^{(k)}$ and
$\theta_N = \theta_N^{(k)} +  
\int_{0}^{1} \e ^{{\bigg[}\int_{s}^{1}A_k(v)\,\eud v{\bigg]}} \cdot 
   \, \vect{b}_k(s) \, \eud s$, and 
 $h_k \geq 0$ is a small multiplier 
that can be found, for example, by
binary search method in the interval
$[0, \, M]$ for some constant $M > 0$.
However, it is important to constrain
$h_k$ to be close to $0$, since
$\eta_k$ must be restricted such that
$\vect{y}_k(t)$ never leaves a tube-like
set that can be determined for
convergence of the Newton-Raphson
method. Alternately, it is also possible
to find $\eta^{\ast}$ such that
$\nabla_{\eta} \beta^2 = 0$, where 
$\beta$ is as in (\ref{Beta-k+1-B}).
In this case, the update in $\eta_k$
is found by the following iteration:
$\eta_0 = 0$ and for $k =1, \, 2, \, \ldots$,
\begin{equation}
\eta_{k+1} = \eta_k -\beta(\theta_0, \, \theta_N)
H^{-1}(\eta_k)  {\bigg[}
\partial_{\theta_0} \beta(\theta_0,\, \theta_N)+
\e ^{{\bigg[}\int_{0}^{1}A_k(v)\,\eud v{\bigg]}}\cdot 
\partial_{\theta_N} \beta(\theta_0,\, \theta_N){\bigg]}
\end{equation}
where $H(\eta)$ is the Hessian matrix of
$\beta^2$ with $\beta$ as in (\ref{Beta-k+1-B}),
$\theta_0 = \theta_0^{(k)}$ and
$\theta_N = \theta_N^{(k)}$.
The method is fast and converges to the
true boundary values provided the initial
approximation is sufficiently accurate.

\section{Summary of Boundary Value Problem and its Extension to Optimal Control Problem and Multidimensional Cases}
In this section, we summarize the dynamic programming method, with an error function using uniform metric. The dynamic programming method can be extended to a vary large class of metrics. We bring out the essential characteristic that is needed for formulation of a dynamic programming based discrete optimization method.  A special formulation for the optimal control problem is given. 

Towards the end of the section, we briefly indicate how to extend the dynamic programming method for the partial differential equations defined in the interior of a simply connected compact, preferably convex, domain, with a regular boundary together with prescribed conditions on it that a solution must satisfy. It is assumed that the boundary is prescribed by a covering as in an atlas. Applications of this method for solving partial differential equations can be found in remote sensing, spectroscopy and tomography, in which regions of physical matter of different permeating, penetrating, reflexivity or resistivity properties, that can affect a flow field, are estimated by the modeling parameters, using measurements taken at the surface and some interior points, where the measurements at these interior points may be assumed or default values, and comparing the reconstruction with another model, which is assumed to be free from anomalies.

\subsection{Summary of Dynamic Programming Method for the Solution to a Boundary Value Problem}

The method described in Section \ref{Section-DPM} can be applied with any nonnegative error functional
$\Error$ for which the $k$-step error functional $\Error_k$ can be evaluated based on $\Error_{k-1}$
and $\error_k$. In particular, if for some function $g$ defined on $[0,\, \infty)\times[0,\, \infty)$,
the following holds
\begin{equation}
\Error_k(\pi_0,\, \ldots,\,\pi_{k-2},\, \pi_{k-1}) =
 g(\Error_{k-1}(\pi_0,\, \ldots,\,\pi_{k-2}),\, \error_k(\pi_{k-1}))
\label{General-k-Step-Error-Recurrence}
\end{equation}
and $\Error_0$ is specified,  then the
dynamic programming method can be
still applied with the function $g$ in stead of $+$.

One of the important error functionals is with
respect to uniform metric, which takes the form
\begin{equation}
\Error = \max_{\scriptstyle{^{t_{i-1} \leq t \leq t_i}_{~ 1 \leq i \leq N}} }
\{ \|A_{i-1}(t) \vect{y}(t) + \vect{b}_{i-1}(t) - \field{f}(\vect{y}(t),\, t)\| \}
\label{Error-WRT-UniformMetric}
\end{equation}
The $k$-step error funcitonal for
(\ref{Error-WRT-UniformMetric})
is given by
\begin{eqnarray}
\error_k(\pi_{k-1}) & =&  \max_{t_{k-1} \leq t \leq t_k}
\|A_{k-1}(t) \vect{y}(t) + \vect{b}_{k-1}(t) - \field{f}(\vect{y}(t),\, t)\| \nonumber \\
\Error_k(\pi_0,\, \ldots\, \pi_{k-2},\, \pi_{k-1}) &  = &
 \max\{\Error_{k-1}(\pi_0,\, \ldots\, \pi_{k-2}),\, \error_k(\pi_{k-1})\}
 \label{k-Step-Error-WRT-UniformMetric}
\end{eqnarray}
where $1 \leq k \leq N$. The error functional (\ref{Error-WRT-UniformMetric})
is best suited for obtaining an initial approximation, followed by the 
Newton-Raphson method of iterative improvement. A bound on the error in the initial
approximation for convergence of the Newton-Raphson mothod for computing the solution of algebraic
or analytic equations can be explicitly found. Now as the formulation for the differential equation that
the error correction term $\vect{y}_k$ satisfies in Sections \ref{Section-Improvements}
and \ref{Section-Modified-Improvements} is obtained by treating the variable $t$ fixed, the same
bound works for convergence of $\vect{y}_k$ (or $\vect{y}_k'$).

Further, with the error functional given by (\ref{Error-WRT-UniformMetric}) the dynamic
programming method can be applied repeatedly until the desired precision is achived using
finer and finer discretizations and restricting the search to only the tube-like set around the
previously obtained optimal parameters. Decisions concerning which parameters (not necessarily
optimal with respect to $S_k$) to retain in subsequent iterations can be made based on
quantitative measures such as stiffness at the parameter value for identifying
the tube set locally. A measure similar to stiffness is the spectrum of
the matrix $A_k(t)$ defined in (\ref{TaylorSeriesApproximation-f}),
which indicates as the iteration progresses which components move
inward and which outward of the tube. The spectrum of  $A_k'(t)$ indicates
the torsion and oscillatory properties of the solution. These measures can
be evaluated at the parameter value from the vector field $\field{f}$ without requiring
complete or even part of the actual solution. Thus besides the tables for $S_k$,
auxiliary tables containing information regarding stiffness or other measures can be
used for choosing the parameter values in the dynamic programming algorithm.
The usefulness of the auxiliary tables is especially significant when using
larger discretization step-sizes. The auxiliary information can expedite
the search method by eliminating unwanted parameter values and 
retaining only those parameter values that could actually produce
the true solution, so that in the subsequent iterations,  closer
approximations (i.e., with smaller error) are produced.
The basic abstract model for the error functional in (\ref{General-k-Step-Error-Recurrence})
can be recast in such a way that allowances for parameter dependences in the accumulating function $g$
and dependence of the $k$-step cost function $\Error_k$ on a future state to reach,
which is described by the model parameter $\pi_{k}$,  for  $1 \leq k  \leq N$, 
are explicated,  as follows:
\begin{equation}
\Error_k(\pi_0,\, \ldots,\,\pi_{k-1},\, \pi_{k}) =
 g(\Error_{k-1}(\pi_0,\, \ldots,\,\pi_{k-1}),\, \error_k(\pi_{k-1}, \pi_{k}))
\label{General-k-Step-Error-Recurrence-current-future}
\end{equation}
and $\Error_0(\pi_0)$ is specified. The model parameter sequence
$(\pi_0,\, \ldots,\,\pi_{k-1},\, \pi_{k})$  in the arguments of the $g$-function in 
(\ref{General-k-Step-Error-Recurrence-current-future}) allows the designer to take into
consideration costs incurred due to lag or drag involved in the course from an initial or past state
up to until the current state, described by the model parameter $\pi_{k-1}$, to reach the next state
described by the model parameter $\pi_{k}$, for $1 \leq k  \leq N$. 

\subsection{Extension to the Solution of Optimal Control Problem}

In this section, we describe a method to solve the optimal control problem
\begin{equation}
\left.
\begin{array}{cl}
\textsf{minimize} ~~ {\mathcal E}(\vect{u}) ~~ \textsf{subject to}\\
   \vect{x}'(t) = \field{f}(\vect{x}(t),\, t) + \vect{u}(t)\,, & \textrm{ for }   0 < t < 1,\\
   \beta(\vect{x}(0), \, \vect{x}(1)) = 0  &  
\end{array}
\right\} \label{Optimal-Control-Problem}
\end{equation}
where $\field{f} \in \Continuous(\reals^n \times [0,\, 1]; \, \reals^n)$,
$\beta \in \Continuous(\reals^n \times \reals^n ;\,  \reals)$,  and
${\mathcal E}(\vect{u}) $ is a cost functional that can be written as
sum of cost functionals over the intervals $[t_{i},\, t_{i+1}]$, $0\leq i \leq N-1$.
Taking $\vect{u}(t) = A_{i}(t)\vect{x}(t)+b_{i}(t)-\field{f}(\vect{x}(t),\, t) $, on
$[t_{i},\, t_{i+1}]$, $0\leq i \leq N-1$,  the dynamic programming method and
improvements of approximations described in the previous section  can be
utilized to find an approximate optimal solution for the control function $\vect{u}(t)$.
For other applications, for instance, in \cite{Varaiya1998}, the optimal control based formulation
is used for estimation of sets reachable from an initial state.

\subsection{Extension for Boundary Value Problems for Partial Differential Equations Defined over a Simply Connected Compact Domain of the Euclidean Space}
In this subsection, we describe a method to solve a boundary value problem of a partial differential equation, defined over a simply connected compact domain, which is a subset of $\reals^{d}$, for some positive integer $d \geq 2$. For simplicity of description, the domain is further assumed to be convex, although this assumption is not always necessary.  The boundary is assumed to be sufficiently regular and specified by smooth surfaces, parameterized as in an atlas, such as in the case of a sphere. The boundary is then propagated inward by, for example, computing the Euclidean distance from a point inside the domain. This approach is called boundary propagation or front propagation \cite{MSV1995}. A change of coordinates, that is consistent with the description of the domain as consisting of concentric surfaces diffeomorphic to the initial boundary of the domain, is performed in the given partial differential equation. These surfaces are parametrically described in such a way that some geometric attribute remains constant for each surface, and hence they are called level sets of the propagating boundary. If the initial boundary is convex, with its unit normal pointing outward, then the boundary can be propagated inward by collecting points obtained by subtracting vectors normal to the boundary at the initial points on the boundary with small multipliers as absolute value from those points on  the initial boundary, as in the gradient descent method, with varying multipliers.  But this propagation may not eventually end in a single point, and, depending on the absolute values chosen as multipliers of the surface normal vectors, various shapes can be realized. In order to construct a set of concentric surfaces ending in a point, the Euclidean distance is computed from a point, which is termed as the center of the domain, to various points on the boundary, whose normal always points outwards, by assumption.  A criterion for an interior point to be the center is stated as follows: let the distance of an interior point to the boundary be defined to be the maximum distance from it to any point on the boundary, and let the center be chosen to be the interior point that results in a minimum distance to the boundary among all interior points. This method of choosing the center is appropriate with compact convex regular boundaries with simply connected interior domains. Then, from the center, the distance to a point on the boundary is multiplied by a scale, as in the case of projective coordinate system, but the scale parameter is chosen to be a real number between 0 and 1. The level sets are the surfaces corresponding to the same scale. Now, appropriate dynamic programming tables are constructed that approximate a solution and its partial derivatives with respect to the changed coordinate system. In this case, the objective is to find a solution that agrees at the center of the domain, when approached from various directions, with minimum error. 

\section{Conclusions}
In this paper, we have presented a dynamic programming based 
formulation for obtaining a best piecewise linear approximation
for continuous vector fields with the solution constrained by
arbitrary boundary conditions. The method attempts obtain a 
fast but reasonably accurate approximate solution. The
method also assumes a discrete space. The parameters to be
determined are of the form  $\pi_i = (A_i(t), \,\vect{b}_i(t),\, \theta_i)$ 
used to approximate the vector field in the interval  $[t_i,\, t_{i+1}]$, for
$0\leq i \leq N-1$. The components $A_i(t)$ and $\vect{b}_i(t)$ contain
undetermined linear combinations of some fixed basis functions (such
as polynomials) defined on the interval $[t_i, \, t_{i+1}]$. The objective function
for minimization of error can be chosen to be the standard $\SquareIntegrable$-norm
or it can also be chosen from among a large class of error functionas including
the uniform metric. The particular aspect of the error functional that allows
a formulation of the dynamic programming method is usually a recurrence relation
that an associated $k$-step functional satisfies. In the proposed method,
the $k$-step functionals satisfy a one-step recurrence relation.
It is possible to formulate dynamic programming based methods also
for error functionals satisfying more general recurrence relations, involving
difficult parameter dependencies. A method for improving an initial approximation
by successive iterations is also presented. The proposed method also updates 
the boundary values. The correction function in the interior of the interval is found by
Newton-Raphson method. The instant when to switch from the dynamic programming
method to successive iterations method for improvement of the current solution can be
determined based on the width of the convergence set for the Newton-Raphson method.
The instant when to switch from the dynamic programming method to successive iterations method for improvement of the current solution can be determined based on the width of the convergence set for the Newton-Raphson method, when it is preferred instead of gradient descent method, for this purpose. It is also possible to consider taking a convex combination of formulated updations, for choosing the actual update.  If any one of the updation formula always produces a more accurate solution than the other, then the convex combination degenerates into binary exclusive combination.


\begin{thebibliography}{99}


\bibitem{RK1965} 
Bellman, Richard, and Robert E. Kalaba, {\em Dynamic programming and modern control theory},
Academic Press, New York, 1965, pp. 195--209


\bibitem{Daniel1978} J. W. Daniel, {\em A road map of methods for approximating solutions of
two point boundary value problems},
in ``Codes for boundary value problems in ordinary differential equations'', 
Bart Childs, Melvin R. Scott, James W. Daniel, Eugene D. Denman, Paul Nelson (Eds.),
Proceedings of a Working Conference May 14-17, 1978, LNCS Vol. 76, Springer Verlag, 1979


\bibitem{Finlayson1972} Bruce A. Finlayson,
{\em The method of weighted residual and variational principals},
Academic Press, 1972


\bibitem{Fox1990} Leslie Fox, {\em The numerical solution of two-point boundary problems
in ordinary differential equations}. Dover Publications, 1990


\bibitem{GR1996}	C. M. Garcia-L\'opez, and J. I. Ramos,
{\em A piecewise-linearized method for ordinary differential equations: two-point boundary-value problems},
 International Journal of Numerical Methods in Fluids, vol. 22, Issue 11, pp.1089-1102

 
\bibitem{SUNDIALS}	Alan.C. Hindmarsh, Peter N. Brown, Keith E. Grant, Steven L. Lee, Radu Serban, Dan E. Shumaker and Carol S. Woodward, {\em SUNDIALS: Suite of nonlinear and differential/algebraic equation solvers}, ACM Transactions on Mathematical Software, Special issue on Advanced CompuTational Software (ACTS),
Vol. 31, Issue 3, September 2005


\bibitem{Keller1987}  Herbert B. Keller, {\em Numerical solution of two point boundary value problems}
Vol. 24. Society for Industrial Mathematics, 1987


\bibitem{MSV1995} 	R. Malladi, J. A. Sethian  and B. C, Vemuri,  {\em Shape modeling with front propagation: a level set approach}, IEEE Trans. on PAMI,  vol 17(2), 1995, pp. 158-175


\bibitem{NS1972}
C. P. Neuman and A. Sen, {\em  Galerkin's procedure, quasilinearization, and nonlinear boundary-value problems}, Journal of Optimization Theory and Applications, Vol. 9, No. 6 / June 1972, pp. 433-437


\bibitem{Osbourne1969}	Osborne, M.R., 
 {\em On Shooting methods for boundary value problems},
  – J. Math. Anal. Appl. 27, 417-433 (1969)


\bibitem{RR2001}	L. S. Ramachandra and D. Roy, {\em A new method for nonlinear
two-point boundary value problems in solid mechanics},
Journal of Applied Mechanics, September 2001, Volume 68, Issue 5, pp. 776-786


\bibitem{Ramos2004} J. I. Ramos, 
{\em Piecewise quasilinearization techniques for singular boundary-value problems}
Computer Physics Commmunicaions, Vol. 158, Issue 1, March 2004, pp. 12-25


\bibitem{RS1972} Sanford M Roberts, and Jerome S. Shipman,
{\em Two-point boundary value problems: shooting methods}
(Modern analytic and computational methods in science and mathematics),
American Elsevier Pub. Co., First Edition edition, 1972


\bibitem{Varaiya1998}	P. Varaiya,
{\em Reach set computation using optimal control},
Proceedings of KIT Workshop, Verimag, Grenoble, 1998.


\end{thebibliography}
\end{document}